\documentclass{scrartcl}

\usepackage[english]{babel}
\usepackage[utf8]{inputenc}
\usepackage{amsmath}
\usepackage{amssymb}
\usepackage{bbm}
\usepackage{mathtools}
\usepackage{xcolor}
\usepackage{graphicx}
\usepackage{ntheorem}

\newcommand{\OFP}{(\Omega, \mathcal{F}, \mathbb{P})}
\newcommand{\NNN}{\mathbb{N}}
\newcommand{\RRR}{\mathbb{R}}
\newcommand{\cL}{{\cal L}}
\newcommand{\cM}{{\cal M}}

\newcommand{\cR}{{\cal R}}
\newcommand{\eins}{\mathbbm{1}}
\newcommand{\ex}{\mathbb{E}}
\newtheorem{remark}{Remark}
\newtheorem{Def}{Definition}
\newtheorem{Prop}{Proposition}
\newtheorem{corollary}{Corollary}
\newtheorem{example}{Example}
\newtheorem{assumption}{Assumption}
\newtheorem{growth condition}{growth condition}
\newtheorem*{proof}{Proof}
\newtheorem{lemma}{Lemma}
\newtheorem{theorem}{Theorem}
\newtheorem{proposition}{Proposition}

\title{Weak continuity of risk functionals with applications to stochastic programming}
\author{M. Claus \footnotemark[2]\ , V. Kr\"atschmer \footnotemark[2] \and R. Schultz \footnotemark[2]}
\begin{document}

\maketitle

\renewcommand{\thefootnote}{\fnsymbol{footnote}}
\footnotetext[2]{Faculty of Mathematics, University of Duisburg-Essen, Campus Essen, Thea-Leymann-Str. 9,
D-45127 Essen, Germany, [matthias.claus][volker.kraetschmer][ruediger.schultz]@uni-due.de}
\renewcommand{\thefootnote}{\arabic{footnote}}

\begin{abstract}
Measuring and managing risk has become crucial in modern decision making under stochastic uncertainty. In two-stage stochastic programming, mean risk models are essentially defined by a parametric recourse problem and a quantification of risk. From the perspective of qualitative robustness theory, we discuss sufficient conditions for continuity of the resulting objective functions with respect to perturbation of the underlying probability measure. Our approach covers a fairly comprehensive class of both stochastic-programming related risk measures and relevant recourse models. Not only this unifies previous approaches but also extends known stability results for two-stage stochastic programs to models with mixed-integer quadratic recourse and mixed-integer convex recourse, respectively.
\\
\medskip
\\
\textit{AMS:} 90C15, 91B30
\end{abstract}

\section{Introduction}

Since the last decade risk management has become an important issue from a practical view point, and as a research field as well.  The interests range from pragmatic solutions for practitioners to research which has founded a hybrid of a new mathematical discipline integrating several fields such as stochastics (e.g. \cite{McNeilEmbrechtsFrey2005}, \cite{Rueschendorf2013}), optimization (e.g. \cite{PflugRoemisch2007}, \cite{Shapiro-etal2009}), numerical analysis. 
(e.g. \cite{Conti-etal2011}) and, when integer variables occur, also algebra and discrete mathematics (e.g. \cite{Schultz2003b}). Since economic risks like credits, prices of stocks or insurance claims are typically faced with uncertainty, most of the methods are settled within a stochastic framework representing risks in terms of random variables. Then basic objects are often stochastic functionals, i.e. real-valued functions defined on sets of random variables expressing economic risks. As a prominent example the so called coherent risk measures may be pointed out. This concept was introduced in  \cite{Artzneretal1999} as a mathematical tool to assess the risks of financial positions.  They are building blocks in quantitative risk management (see \cite{McNeilEmbrechtsFrey2005}, \cite{PflugRoemisch2007}, \cite{Rueschendorf2013}), and they have been suggested as a systematic approach for calculations of insurance premia
(cf. \cite{Kaasetal2008}). Besides the ordinary expectation, the conditional value at risk and the upper semideviation are the most known examples for coherent risk measures. However, meanwhile the more general notion of convex risk measure has replaced coherent risk measures in playing their roles.

\medskip

Of particular interest are stochastic functionals which are distribution invariant,  identifying risks with identical distributions. For instance the expectation, the conditional value at risk and the upper semideviation satisfy this property. They all may be redefined as functionals on sets of probability measures representing the distributions of the risks. Recent contributions analyze analytic properties of such functionals like specific types of continuity and differentiability (cf. \cite{KraetschmerEtAll2014}, \cite{KraetschmerEtAll2015}). Such properties have immediate applications for statistical issues of the functionals e.g. the sample average approximation method (SAA) (\cite{KraetschmerEtAll2014}, \cite{KraetschmerEtAll2015} again; see also \cite{PflugRoemisch2007}, \cite{Shapiro-etal2009}, \cite{BelomestnyKraetschmer2012}, \cite{Rueschendorf2013}). The aim of the present paper is to point out how investigations of stochastic programming problems may profit from continuity properties of distribution invariant stochastic functionals which are used for objective functions. In technical terms, investigations are devoted to topological considerations on spaces of random variables induced in a natural way by stochastic programming problems. 
 
\medskip

Stochastic programming is based on the crucial assumption that uncertainty can be captured by a probability measure which, in turn, 
has impact on structural and/or algorithmic properties of  the objective function and/or the constraints. The probability measure usually being subjective or resulting from statistical estimation the issue of stability comes to the fore, i.e., small perturbations of the measure shall lead to only small perturbations of the optimal value and the optimal solution sets. Beginning with \cite{BirgeWets1986}, one line of research in stochastic programming has addressed questions of stability in the theoretical framework of epi-convergence (cf. e.g. \cite{Attouch1984}, \cite{AttouchWets1991}, \cite{RockafellarWets1998}). However, for the problems considered in this paper, a different approach based on arguments from \cite{Berge1963} appears to be more straightforward (see Remark \ref{RemEpi} after Corollary \ref{CorStability} for details). It appeals to parametric optimization by aiming at (semi-)continuity of optimal-value functions and of multifunctions given by sets of optimal solutions. Typically, the parameter spaces may vary from Euclidean spaces of parameters of distributions to topological or metric spaces of probability measures equipped with weak convergence of probability measures \cite{Kall1987,RobinsonWets1987} or with suitable probability metrics \cite{RoemischWakolbinger1987, RachevRoemisch2002, Roemisch2003}. In classical stochastic programming the objective function is described in terms of expectations representing a risk neutral attitude of the decision maker. More recent contributions try to incorporate risk aversion of the decision makers using different distribution invariant coherent risk measures, where also both, continuous or mixed-integer variables are involved (cf. e.g.\cite{RoemischVigerske2008,MaerkertSchultz2005,SchultzTiedemann2006,{Shapiro-etal2009}}).

\medskip

In the nutshell stability of stochastic programming refers to continuity of the distribution invariant functionals used for the objective functions. In the above mentioned literature this viewpoint has not been exploited systematically. Instead by individual reasoning suitable settings of uniform integrability  or moment conditions had to be adapted to the individual objectives (starting with \cite{Kall1987,RobinsonWets1987,RoemischWakolbinger1987}), and sufficient conditions had to be found to transfer weak convergence of sequences of probability measures to sequences of image measures (see \cite{Schultz2003a} for an example of detailed elaboration).

\medskip

The present paper is an attempt of systemization. We provide an umbrella for most of the different settings in the  papers referred to in the previous paragraphs. The line of reasoning is inspired by recent investigations on continuity of distribution invariant convex risk measures (\cite{KraetschmerEtAll2014}, \cite{KraetschmerEtAll2015b}). They are based on the so called $\psi-$weak topologies which is a quite new class of topologies for sets of probability measures enclosing the topology of weak convergence. We shall extend the studies to more general distribution invariant stochastic functionals which will be referred to as risk functionals.

\medskip 

More precisely, we introduce general risk functionals living on sets of probability measures satisfying some moment condition and allowing for specification by proper choices of integrands. Every such moment condition corresponds with some particular $\psi-$weak topology, and as a basic observation, all the considered risk functionals are continuous w.r.t. the related $\psi-$weak topologies. We shall then identify sufficient growth conditions to integrands of risk functionals implying along with the continuity of the risk functionals the continuity of resulting  objective function with respect to some particular $\psi-$weak topology induced by the growth condition. In general, $\psi-$weak topologies are finer than the topology of weak convergence. However, we may specify exactly those subsets, where they coincide. Hence in the last step we may point out the domains of stability for stochastic programming involving the considered general objective functions.

\medskip

We shall apply the technical results to two-stage mean-risk models unifying previous work. Moreover, this top-down approach  foremostly yields verifiable continuity conditions for broader classes of risk functionals than before. For stochastic programs this enables extension of the continuity, and thus stability, analysis to more comprehensive classes of models.

\medskip

The paper is organized as follows. In section 2, we provide a unifying view on a class of stochastic optimization problems enclosing the two-stage case and various notions of risk aversion. Our main result is on weak continuity of the objective function w.r.t. the decision and the underlying probability measure and allows for conclusions about qualitative stability. The approach is applicable to two-stage problems whenever the optimal value function of the recourse problem is Borel measurable, polynomially bounded in terms of the entering parameters and continuous outside of a suitable set. In section 3, we check these conditions for various recourse models: For stochastic programs with linear or mixed-integer linear recourse, the conditions hold under standard assumptions and allow to unify various existing proofs of stability. Furthermore, we extend the analysis to the cases of mixed-integer quadratic and mixed-integer convex recourse. Section 4 is devoted to $\psi-$weak topologies that are an important tool in our argumentation. Their relationship with the topology of weak convergence will also be discussed there. Finally, we shall be ready to proof the main result in Section 5.

\section{Main result} \label{mainresult}

Let $X \subseteq \mathbb{R}^{n}$ be a nonempty set, $f: \mathbb{R}^{n} \times \mathbb{R}^{s} \to \mathbb{R}$ a Borel-measurable mapping and let $Z: \Omega' \to \mathbb{R}^{s}$ denote a fixed and known random vector on some probability space $(\Omega', \mathcal{F}', \mathbb{P}')$. We consider the stochastic programming problem
\begin{equation}
\label{OptimizationUnderUncertainty}
\min_{x} \lbrace f(x,Z(\omega)) \; | \; x \in X \rbrace,
\end{equation}
where the decision on $x$ has to be made nonanticipatorily of the realization $Z(\omega)$. Two-stage problems arise from \eqref{OptimizationUnderUncertainty} if $f$ itself is given by the optimal value function of an optimization problem parametrized in $x$ and $Z(\omega)$.

\medskip

While not well defined because of the nonanticipativity constraint, \eqref{OptimizationUnderUncertainty} may be understood as selecting in some sense a minimal random variable out of the family
\begin{equation*}
f(X,Z) := \lbrace f(x,Z(\cdot)) \; | \; x \in X \rbrace.
\end{equation*}
The notion of minimality can be specified by endowing $f(X,Z)$ with a preorder, i.e. a reflexive and transitive, yet not necessarily antisymmetric binary relation. Provided that the members of $f(X,Z)$ are all integrable, we might rank them by their expectations. This approach leads to the risk neutral model
\begin{equation}
\label{RiskNeutralProblem}
\min_{x} \lbrace \mathbb{E}[f(x,Z)] \; | \; x \in X \rbrace.
\end{equation}
In many applications, the random vector $Z$ may be subject to perturbations. This motivates stability analysis for the problem's optimal value and its set of optimal solutions. Since the only relevant information of $Z$ for problem \eqref{RiskNeutralProblem} is provided by its distribution, we might equivalently work with perturbations of this distribution in the space of Borel probability measures on $\RRR^{s}$. For qualitative stability analysis, we endow this space with the topology of weak convergence (see e.g. \cite{Billingsley1968}).

\medskip

The stochastic programming problem \eqref{RiskNeutralProblem} then reads
\begin{equation}
\label{RiskNeutralDistribution}
\min_{x} ~\left\lbrace \int_{\RRR^n \times \RRR^s} f~d~\delta_{x}\otimes\nu_{Z} \; \big| \; x \in X \right\rbrace,
\end{equation}
where $\delta_{x}\otimes\nu_{Z}$ denotes the product measure of the Dirac measure at $x$ and the distribution $\nu_{Z}$ of $Z$. Qualitative stability of this stochastic programming problem at $\nu_{Z}$ is entailed by continuity of 
\begin{equation*}
(x,\nu) \mapsto \int_{\RRR^n \times \RRR^s} f~d~\delta_{x}\otimes\nu = \int_{\RRR} t~(\delta_{x}\otimes\nu)\circ f^{-1}(dt)
\end{equation*}
at $(x,\nu_{Z})$ for every $x\in X$ w.r.t. the product topology of the standard topology on $\RRR^{n}$ and the topology of weak convergence on the space of all $s$-dimensional distributions (see Corollary \ref{CorStability} below).

\medskip

If we want to take into account risk aversion we might change the objective function in \eqref{RiskNeutralDistribution} to
$$
\cR\big((\delta_{x}\otimes\nu_{Z})\circ f^{-1}\big),
$$
where $\cR$ denotes a functional on a set ${\cal N}$ of Borel probability measures on $\RRR$. Typically, functionals $\cR$ which are nondecreasing w.r.t. the increasing convex order may be considered as suitable choices. Outstanding examples are provided by functionals which are derived from so called law-invariant convex risk measures. More precisely, let $(\Omega, \mathcal{F}, \mathbb{P})$ be an atomless probability space, i.e. it supports some random variable $U$ which is uniformly distributed on $]0,1[.$ Furthermore, for $p\in [1,\infty[$ denote by $L^{p}(\Omega, \mathcal{F}, \mathbb{P})$ the standard $L^{p}-$space w.r.t. $\OFP$. Then a mapping $\rho: L^{p}(\Omega, \mathcal{F}, \mathbb{P}) \to \mathbb{R}$ is called a \textit{convex risk measure} if it is a convex function which is nondecreasing w.r.t. the $\mathbb{P}$-almost sure partial order and \textit{translation-equivariant}, i.e. $\rho(X+t) = \rho(X) + t$ for $X\in L^{p}(\Omega, \mathcal{F}, \mathbb{P})$ and $t\in\RRR$. Popular examples are the following

\medskip

\begin{example}
\label{Average Value at Risk}
Average Value at Risk: $p = 1,$
$$
AV@R_{\alpha}(Y) := \frac{1}{1 - \alpha}~\int_{\alpha}^{1}F_{Y}^{\leftarrow}(\beta)~d\beta\quad\mbox{for some}~\alpha\in ]0,1[,
$$
where $F_{Y}^{\leftarrow}$ denotes the left-continuous quantile function of the distribution function $F_{Y}$ of $Y$.
\end{example}

\medskip

\begin{example}
\label{mean upper semideviation}
Mean upper semideviation of order $p$:
$$
\rho_{a,p}(Y) :=  \ex[Y] + a \left(\ex\left[\Big(\big(Y - \ex[Y]\big)^{+}\Big)^{p}\right]\right)^{1/p}\quad\mbox{for some}~a\in [0,1].
$$
\end{example}

In addition, the convex risk measure $\rho$ is called \textit{law-invariant} if $\rho(Y) = \rho(\widetilde{Y})$ whenever $Y$ and $\widetilde{Y}$ have the same law under $\mathbb{P}.$ Obviously, the Average Value at Risk as well as the mean upper semideviation are law-invariant. Now, any law-invariant risk measure $\rho: L^{p}\OFP\rightarrow\RRR$ is associated with a functional $\cR_{\rho}$ on the set of Borel probability measures on $\RRR$ with absolute moments of order $p$ via 
\begin{equation}
\label{AssociatedFunctional}
\cR_{\rho}(\mu) := \rho(F^{\leftarrow}_{\mu}(U)),
\end{equation}
where $F^{\leftarrow}_{\mu}$ stands for the left-continuous quantile function of the distribution function $F_{\mu}$ of $\mu.$ This functional is nondecreasing w.r.t. the increasing convex order (cf. e.g. Theorem \ref{RRhoStetig} below).

\medskip

It will turn out that for our purposes, we may relax the conditions on the mapping $\rho$ by dropping the requirement of translation-equivariance. So our investigation is built upon the following assumption: 

\medskip

\begin{assumption}
\label{Arho}
$\rho: L^{p}\OFP\rightarrow\RRR$ is
\begin{itemize}
\item a convex function,
\item nondecreasing w.r.t. the $\mathbb{P}$-almost sure partial order,
\item law-invariant,
\item defined on the $L^{p}-$space of the atomless probability space $\OFP$ and $p \geq 1$.
\end{itemize}
\end{assumption}

\medskip

Let us mention an example.

\medskip

\begin{example}
\label{mean upper semideviation target}
Mean upper semideviation of order $p$ from a target:
$$
\rho_{a,c,p}(Y) :=  \ex[Y] + a \left(\ex\left[\Big(\big(Y - c\big)^{+}\Big)^{p}\right]\right)^{1/p}\quad\mbox{for some}~a\in [0,1],~c > 0,
$$
which does not fulfill translation-equivariance (cf. \cite[Example 6.25]{Shapiro-etal2009}).
\end{example}

\medskip

We associate the law-invariant mapping $\rho$ with the functional ${\cal R}_{\rho}$ defined in \eqref{AssociatedFunctional} and consider the parametric stochastic programming problem 
\begin{equation}
\label{StochasticProblem}
\min_{x}\{Q(x,\nu)\mid x\in X\},
\end{equation}
where 
\begin{equation}
\label{DefQ}
Q(x,\nu) := {\cal R}_{\rho}\left((\delta_{x}\otimes\nu)\circ f^{-1}\right).
\end{equation}
$Q$ is well defined provided that the image probability measure of $\delta_{x}\otimes\nu$ under $f$ has absolute moments of order $p$. The latter has to be guaranteed by imposing an additional assumption on the growth of $f$.

\medskip

\begin{Def}
We call a mapping $e: \RRR^n \times \RRR^s \to \RRR$ limited by the exponent $\gamma \in ]0,\infty[$ iff there is some locally bounded mapping $\eta:\RRR^{n} \to ]0,\infty[$ such that
$$
|e(x,z)| \leq \eta(x)(\|z\|^{\gamma} + 1)\quad\mbox{for all}~(x,z)\in\RRR^{n}\times\RRR^{s},
$$
where $\|\cdot\|$ denotes the Euclidean norm on $\mathbb{R}^{s}$.
\end{Def}

\medskip

We shall work with the following assumption:

\medskip

\begin{assumption}
\label{Af}
$f: \RRR^n \times \RRR^s \to \RRR$ is
\begin{itemize}
\item Borel-measurable,
\item limited by an exponent $\gamma \in ]0,\infty[$. 
\end{itemize}
\end{assumption}

\medskip

Denoting by ${\cal N}$ some suitable subset of Borel probability measures on $\RRR$ specified below, it will be shown later on that Assumptions \ref{Arho} and \ref{Af} imply the continuity of the mapping
$$
Q:\RRR^{n}\times \mathcal{N} \rightarrow \mathbb{R},~Q(x,\nu) = {\cal R}_{\rho}\left((\delta_{x}\otimes\nu)\circ f^{-1}\right)
$$ 
w.r.t. the product topology of the standard topology on $\RRR^{n}$ and some particular topology on ${\cal N}$ which is finer than the ordinary topology of weak convergence. It is a special case of the so called $\psi$-weak topology which we shall recall next.

\medskip

For this purpose let $\psi:\RRR^{d}\to [0,\infty[$ be a continuous function such that $\psi\geq 1$ outside a compact set with some $d\in \mathbb{N}$. Such a function will be referred to as a {\em gauge function}. Let ${\cal M}_1^\psi(\RRR^{d})$ be the set of all Borel-probability measures $\mu$ on $\RRR^{d}$ satisfying $\int\psi\,d\mu<\infty$, and $C_\psi(\RRR^{d})$ be the space of all continuous functions $g:\RRR^s \rightarrow \RRR$ for which $\sup_{y\in\RRR^s}|g(y)/(1+\psi(y))|<\infty$. The $\psi$-weak topology on ${\cal M}_1^\psi(\RRR^{d})$ is defined to be the coarsest topology for which all mappings $\mu\mapsto\int g\,d\mu$, $g\in C_\psi(\RRR^{d})$, are continuous; cf.\ Section A.6 in \cite{FoellmerSchied2011}. Since $C_\psi(\RRR^{d})$ encloses all bounded continuous real-valued mappings on $\RRR^{d}$, the $\psi$-weak topology is finer than the topology of weak convergence on $\cM_{1}^{\psi}(\RRR^d)$.

\medskip

Choosing $d = 1$ and $\psi = |\cdot|^{p}$, we may describe formally the functional $\cR_{\rho}$  by
$$
\cR_{\rho}:\cM_{1}^{|\cdot|^{p}}(\RRR)\rightarrow\RRR,~\cR_{\rho}(\mu) = \rho(F^{\leftarrow}_{\mu}(U)).
$$

Moreover, under Assumptions \ref{Arho} and \ref{Af}, we may choose ${\cal N} = \cM_{1}^{\|\cdot\|^{\gamma p}}(\RRR^{s})$ to define the domain of the function $Q$.

\medskip

\begin{lemma} \label{LemmaEndlichkeit}
Under Assumption \ref{Af}, $\int_{\mathbb{R}^{s}} |f(x,z)|^{p}~\nu(dz) < \infty$ holds for any $x \in \mathbb{R}^{n}$ and any $\nu \in \mathcal{M}^{\| \cdot \|^{\gamma p}}_{1}(\mathbb{R}^{s})$.
\end{lemma}

\medskip

\begin{proof}
By Assumption \ref{Af} there exists a locally bounded mapping $\eta: \mathbb{R}^{n} \to ]0,\infty[$ such that
$$
|f(x,z)|^{p}\leq \eta(x)^{p}(\|z\|^{\gamma} + 1)^{p}\leq \eta(x)^{p} 2^{p}\max\{\|z\|^{p\gamma},1\}\leq \eta(x)^{p} 2^{p}(\|z\|^{p\gamma} + 1)
$$
holds for any $(x,z) \in \RRR^{n} \times \RRR^s$. The statement of Lemma \ref{LemmaEndlichkeit} follows immediately.
\end{proof}

\medskip

The line of our reasoning will be as follows. First, we shall prove in Theorem \ref{QStetig} below that under Assumptions \ref{Arho} and \ref{Af}, the mapping
$$
Q:\RRR^{n}\times \mathcal{M}^{\| \cdot \|^{\gamma p}}_{1}(\mathbb{R}^{s}) \to \mathbb{R},~Q(x,\nu) = {\cal R}_{\rho}\left((\delta_{x}\otimes\nu)\circ f^{-1}\right)
$$ 
is continuous w.r.t. the product topology of the standard topology on $\RRR^{n}$ and the $\|\cdot\|^{\gamma p}$-weak topology on $\mathcal{M}^{\| \cdot \|^{\gamma p}}_{1}(\mathbb{R}^{s})$. Then the key of our argumentation will be to point out those subsets $\cM$ of $\mathcal{M}^{\| \cdot \|^{\gamma p}}_{1}(\mathbb{R}^{s})$ where the $\|\cdot\|^{\gamma p}$-weak topology coincides with the topology of weak convergence. It will turn out (see Section \ref{psi-weak topology} and in particular Proposition \ref{coincidence}) that exactly the locally uniformly $\|\cdot\|^{\gamma p}-$integrating subsets $\cM$ satisfy this property. Here a subset $\cM\subseteq \mathcal{M}^{\| \cdot \|^{\gamma p}}_{1}(\mathbb{R}^{s})$ will be called \textit{locally uniformly $\|\cdot\|^{\gamma p}-$integrating} if for every $\varepsilon > 0$ and any $\nu\in\cM$ there exists some open neighborhood $\cL$ of $\nu$ w.r.t. the topology of weak convergence such that
$$
\lim_{a\to\infty}\sup_{\mu\in\cL\cap\cM}\int_{\RRR^{s}}~\|z\|^{\gamma p} \cdot\eins_{]a,\infty[}(\|z\|^{\gamma p})~\mu(dz)\leq\varepsilon.
$$ 
This concept has been established recently in \cite{Zaehle2016} as the suitable one to identify sets of Borel probability measures where the $\|\cdot\|^{\gamma p}$-weak topology and the topology of weak convergence are the same. In \cite{KraetschmerEtAll2015b}, further equivalent characterizations for locally uniformly $\|\cdot\|^{\gamma p}$-integrating sets of Borel probability measures have been found, providing also a large amount of examples. Simple examples may be constructed by imposing restrictions on the absolute moments.
\medskip

\begin{example}
\label{moment condition}
Let $\kappa, \varepsilon > 0$. Then by \cite[Corollary A.47, (c)]{FoellmerSchied2011}, the set 
$$
\cM_{1,\kappa}^{\| \cdot \|^{\gamma p + \varepsilon}}(\RRR^{s}) := 
\left\{\mu\in \mathcal{M}^{\| \cdot \|^{\gamma p + \varepsilon}}_{1}(\mathbb{R}^{s})~\big|~ \int_{\RRR^{s}}\| z\|^{\gamma p + \varepsilon}~\mu(dz)\leq\kappa\right\}
$$
is a relatively compact subset of $\mathcal{M}^{\| \cdot \|^{\gamma p}}_{1}(\mathbb{R}^{s})$ for the $\|\cdot\|^{\gamma p}$-weak topology. Then in view of Lemma \ref{relative compactness} below it is also a locally uniformly $\|\cdot\|^{\gamma p}$-integrating subset of $\mathcal{M}^{\| \cdot \|^{\gamma p}}_{1}(\mathbb{R}^{s})$.
\end{example}

\medskip

Now we are ready to formulate our main result. 
\medskip

\begin{theorem} \label{QRestriktionStetig}
Let Assumptions \ref{Arho} and \ref{Af} be fulfilled and $D_f$ denote the set of discontinuity points of $f$. Furthermore, let $\cM\subseteq \mathcal{M}^{\| \cdot \|^{\gamma p}}_{1}(\mathbb{R}^{s})$ be a locally uniformly $\|\cdot\|^{\gamma p}$-integrating subset and $(x,\nu) \in \RRR^n \times \cM$ such that $\delta_{x} \otimes \nu(D_{f}) = 0$. Then the restriction $Q|_{ \mathbb{R}^{n} \times \cM}$ is continuous at $(x,\nu)$ w.r.t. the product topology of the standard topology on $\mathbb{R}^{n}$ and the topology of weak convergence on $\cM$.
\end{theorem}

\medskip

Theorem \ref{QRestriktionStetig} has the following specialization.

\medskip

\begin{corollary} \label{CorLebesgue}
Let Assumptions \ref{Arho} and \ref{Af} be fulfilled and let $\cM\subseteq \mathcal{M}^{\| \cdot \|^{\gamma p}}_{1}(\mathbb{R}^{s})$ be a locally uniformly $\|\cdot\|^{\gamma p}$-integrating subset. Furthermore, let $x \in \mathbb{R}^{n}$ be such that $\{z\in\RRR^{s}\mid (x,z)\in D_{f}\}$ has Lebesgue measure $0$ assume that $\nu \in \cM$ is absolutely continuous w.r.t. the Lebesgue-Borel measure on $\RRR^{s}$. Then $Q|_{ \mathbb{R}^{n} \times \cM}$ is continuous at $(x, \nu)$ w.r.t. the product topology of the standard topology on $\mathbb{R}^{n}$ and the topology of weak convergence on $\cM$.
\end{corollary}

\medskip

\begin{proof}
We have $\delta_{x}\otimes\nu(D_{f}) = 0$, so that Theorem \ref{QRestriktionStetig} is applicable.
\end{proof}

\medskip

Theorem \ref{QRestriktionStetig} has immediate implications towards stability of problem \eqref{StochasticProblem}. Let
$$
\varphi: \mathcal{M}^{\| \cdot \|^{\gamma p}}_{1}(\mathbb{R}^{s}) \to [-\infty,\infty[,~\varphi(\nu) = \inf_{x} \{Q(x,\nu)\mid x\in X\}
$$
and
$$
\Phi: \mathcal{M}^{\| \cdot \|^{\gamma p}}_{1}(\mathbb{R}^{s}) \to 2^{X},~\Phi(\nu) = \{ x \in X \mid Q(x,\nu) = \varphi(\nu) \}
$$
its optimal value function and its optimal solution set mapping.

\medskip

\begin{corollary} \label{CorStability}
Let Assumptions \ref{Arho} and \ref{Af} be fulfilled and let $\cM\subseteq \mathcal{M}^{\| \cdot \|^{\gamma p}}_{1}(\mathbb{R}^{s})$ be a locally uniformly $\|\cdot\|^{\gamma p}$-integrating subset. Furthermore, let $\nu \in \mathcal{M}$ be such that $(\delta_{x} \otimes \nu)(D_{f}) = 0$ holds for any $x \in X$. Then the restriction $\varphi|_{\cM}$ is upper semicontinuous at $\nu$ w.r.t. the topology of weak convergence on $\mathcal{M}$. If $X$ is compact, $\varphi|_{\cM}$ is continuous at $\nu$ and $\Phi|_{\cM}$ is upper semicontinuous at $\nu$ w.r.t. the topology of weak convergence on $\mathcal{M}$. In this case, $\Phi|_{\cM}(\nu)$ is nonempty and compact.
\end{corollary}

\medskip

\begin{proof}
See e.g. \cite[Section 4.1]{BonnansShapiro2000} for the first and \cite[Theorem 2]{Berge1963} for the second part.
\end{proof}

\medskip

\begin{remark} \label{RemEpi}
One could try to obtain a similar result by employing the theory of epi-convergence. Such an approach would only require lower semicontinuity of the functional $Q(\cdot, \mu)$ and might allow for a weaker assumption than $(\delta_{x} \otimes \mu)(D_{f}) = 0$. However, one might have to find additional assumptions under which weak convergence of a sequence $\lbrace \mu_n \rbrace_{n \in \mathbb{N}} \subset \mathcal{M}$ to $\mu \in \mathcal{M}$ entails epi-convergence of the associated sequence of functionals $\lbrace Q(\cdot,\mu_n) \rbrace_{n \in \mathbb{N}}$. In case that $Q(\cdot, \mu)$ is convex, it is sufficient to ensure pointwise convergence of $\lbrace Q(\cdot,\mu_n) \rbrace_{n \in \mathbb{N}}$. Whereas in general, even under the assumptions of Corollary \ref{CorStability}, the situation is less clear. Note that for the recourse models considered in Section 3, $Q(\cdot, \mu)$ may fail to be convex.
\end{remark}

\section{Application to two-stage mean-risk models}

Two-stage stochastic programs can be seen as the special case of \eqref{OptimizationUnderUncertainty} where $f$ is the optimal value function of a recourse problem depending on the parameters $x$ and $Z(\omega)$. While we may use Corollary \ref{CorStability} to derive stability for such problems, the verification of Assumption \ref{Af} becomes an issue. Furthermore, the corollary can only be applied if $(\delta_{x} \otimes \nu)(D_{f}) = 0$ holds for any $x \in X$. Hence, situations in which an explicit description of a suitable superset of the set of discontinuities of $f$ is available are of special interest. In this section we address these issues by providing sufficient conditions for various classes of recourse problems.

\subsection{Linear recourse}

We first turn our attention to the case of a linear recourse problem, i.e. the situation where
\begin{equation}
\label{DeffLin}
f(x,z) = \inf_y \lbrace q(x,z)^\top y \mid Ay = h(x,z), y \geq 0 \rbrace
\end{equation}
for a matrix $A \in \RRR^{k \times m}$ and mappings $q:\RRR^n \times \RRR^s \to \RRR^m$ and $h:\RRR^n \times \RRR^s \to \RRR^k$. Such problems have been featured prominently in the stochastic programming literature and our approach using Corollary \ref{CorStability} will only slightly extend the existing results (see e.g. \cite{RoemischSchultz1993}). Nevertheless, it allows to unify and simplify the proofs of stability for special risk measures (in particular the ones mentioned in Examples \ref{Average Value at Risk} to \ref{mean upper semideviation target}). We shall work with the following assumptions.

\medskip

\begin{assumption} \label{ALin}
\begin{itemize}
\item $A$ has full rank,
\item $\lbrace y \in \RRR^m \mid Ay = h(x,z), y \geq 0 \rbrace \neq \emptyset$ for any $(x,z) \in X \times \RRR^s$,
\item $\lbrace u \in \RRR^k \mid A^{\top} u \leq q(x,z) \rbrace \neq \emptyset$ for any $(x,z) \in X \times \RRR^s$,
\item the mappings $q$ and $h$ are continuous,
\item $\|q(\cdot, \cdot)\|$ and $\|h(\cdot, \cdot)\|$ are limited by exponents $\gamma_q, \gamma_h \in ]0,\infty[$, respectively.
\end{itemize}
\end{assumption}

\medskip

The first three conditions in Assumption \ref{ALin} are standard assumptions in two-stage stochastic programming with linear recourse. Note that the other two conditions are automatically fulfilled if the mappings $q$ and $h$ are linear. 

\medskip

\begin{Prop}
Under Assumption \ref{ALin}, the mapping $f$ defined in \eqref{DeffLin} is finite, continuous and limited by $\gamma_q + \gamma_h$. In particular, Assumption \ref{Af} is fulfilled.
\end{Prop}

\medskip

\begin{proof}
Finiteness and continuity of $f$ are a well known conclusion from the Basis Decomposition Theorem in \cite{WalkupWets1969}. By the same result, there exists a finite number of matrices $B_{1}, \ldots, B_{N} \in \mathbb{R}^{k \times m}$ such that for any $(x,z) \in X \times \RRR^s$, there is an index $i \in \lbrace 1, \ldots, N \rbrace$ satisfying
\begin{equation*}
|f(x,z)| \; = \; |h(x,z)^{\top}B_{i} q(x,z)| \; \leq \; \kappa_{B} \|h(x,z)\| \|\bar{q}(x,z)\|,
\end{equation*}
where $\kappa_{B} := \max_{j = 1, \ldots, N} \| B_j \|_{\mathcal{L}(\mathbb{R}^{m}, \mathbb{R}^{k})} < \infty$ denotes the maximal operator norm among the matrices $B_1, \ldots, B_N$. Invoking the final part of Assumption \ref{ALin} and the fact that
$$
(\|z\|^{\gamma_q} + 1)(\|z\|^{\gamma_h} + 1) \leq 3(\|z\|^{\gamma_q + \gamma_h} + 1)
$$
holds for any $z \in \RRR^s$, we may conclude that $f$ is limited by $\gamma_q + \gamma_h$.
\end{proof}

\medskip

Consequently, Corollary \ref{CorStability} is applicable whenever $f$ is given by \eqref{DeffLin} and Assumptions \ref{Arho} and \ref{ALin} are fulfilled. In this case, $D_f = \emptyset$ and we obtain qualitative stability of problem \eqref{StochasticProblem} on every locally uniformly $\|\cdot\|^{(\gamma_q + \gamma_h) p}$-integrating subset of $\cM_1^{(\gamma_q + \gamma_h) p}$.

\subsection{Mixed-integer linear recourse}

Let us turn over to the case of a mixed-integer linear recourse problem, i.e. let $f$ take the form
\begin{equation}
\label{MILPf}
f(x,z) = \inf_{y} \lbrace  q^{\top}y \mid Ay = h(x,z), y \geq 0, y \in \RRR^{m_1} \times \mathbb{Z}^{m_2} \rbrace
\end{equation}
for a vector $q \in \mathbb{R}^{m_1 + m_2}$, a matrix $A \in \mathbb{R}^{k \times (m_1 + m_2)}$ and a mapping $h: \RRR^n \times \RRR^s \to \RRR$. Note that the vector $q$ in the objective function does not depend on $(x,z)$ anymore. While the analytical properties of the problem in \eqref{MILPf} are weaker than those of the problem in \eqref{DeffLin}, they are still sufficient to derive the conditions of Assumption \ref{Af} under standard assumptions. This allows to apply Corollary \ref{CorStability} to two-stage stochastic programs with mixed-integer linear recourse. Again, such problems have been studied extensively (see e.g. \cite{MaerkertSchultz2005} and \cite{SchultzTiedemann2006}). In the mentioned papers, the analysis is restricted to sets of Borel probability measure having uniformly bounded moments of order striclty greater than $\gamma p$. By Example \ref{moment condition}, such sets are locally uniformly $\|\cdot\|^{\gamma p}$-integrating and hence a special case of the present framework. On top of this generalization, our work provides a unified proof of stability for the special risk measures discussed in Examples 1 to 3. We shall impose the following assumptions.

\medskip

\begin{assumption} \label{AMILP}
\begin{itemize}
\item $A$ has rational entries,
\item $\lbrace y \in \RRR^{m_1} \times \mathbb{Z}^{m_2} \mid Ay = t, y \geq 0 \rbrace \neq \emptyset$ for any $t \in \RRR^k$,
\item $\lbrace u \in \RRR^k \mid A^{\top} u \leq q \rbrace \neq \emptyset$,
\item the mapping $h$ is continuous,
\item $\|h(\cdot, \cdot)\|$ is limited by an exponent $\gamma_h \in ]0,\infty[$.
\end{itemize}
\end{assumption}

\medskip

Note that the second part of Assumption \ref{AMILP} implies that $A$ has full rank.

\medskip

\begin{Prop}
Under Assumption \ref{AMILP}, the mapping $f$ defined in \eqref{MILPf} is finite, upper semicontinuous and limited by the exponent $\gamma_h$. In particular, Assumption \ref{Af} is fulfilled. Let $A_1 \in \mathbb{R}^{k \times m_1}$ be the matrix of the first $m_1$ columns of $A$ and denote by $A_2 \in \mathbb{R}^{k \times m_2}$ the matrix of the last $m_2$ columns of $A$. Then $f$ is continuous outside of
\begin{equation}
\label{DiscMILPf}
h^{-1} \left(\bigcup_{y_2 \in \mathbb{Z}^{m_2}, y_2 \geq 0} (\lbrace A_2 y_2 \rbrace + \mathcal{A}) \right),
\end{equation}
where $\mathcal{A}$ denotes the boundary of the set $\lbrace A_1 y_1 \mid y_1 \geq 0\rbrace$.
\end{Prop}

\medskip

\begin{proof}
Under our standard assumptions, finiteness, upper semicontinuity and statement about the discontinuity points of $f$ follow directly from the classical results in \cite{BankMandel88} and \cite{BlairJeroslow77}. By the same results, there exist constants $\alpha, \beta > 0$ such that
$$
|f(x,z)-f(x',z')| \leq \alpha \|h(x,z) - h(x',z')\| + \beta
$$
holds for any $(x,z), (x',z') \in \RRR^n \times \RRR^s$. From the latter, we may conclude that $f$ is limited by the same exponent as $\|h(\cdot,\cdot)\|$.
\end{proof}

\medskip

The set in \eqref{DiscMILPf} is the preimage of a countable union of hyperplanes under the mapping $h$. In particular, if we have $k=s$ and for any $x \in X$, $h_{x}: \mathbb{R}^s \to \mathbb{R}^s$ is a $C^{1}-$diffeomorphism, the set $\lbrace z \in \RRR^s \mid (x,z) \in D_f \rbrace$ has Lebesgue measure $0$ for any $x \in X$. Consequently, $(\delta_x \otimes \nu)(D_f) = 0$ holds for any $x \in X$ and any Borel probability measure $\mu$ that is absolutely continuous w.r.t. the Lebesgue measure.

\subsection{Mixed-integer quadratic recourse}

In this subsection, we shall consider recourse problems with quadratic objective, linear constraints and mixed-integer variables. Consequently, $f$ takes the form
\begin{equation}
\label{MIQPf}
f(x,z) = \inf_{y} \lbrace  y^{\top}Dy + q(x,z)^{\top}y \mid Ay \leq h(x,z), y \in \RRR^{m_1} \times \mathbb{Z}^{m_2} \rbrace,
\end{equation}
where $q: \RRR^n \times \RRR^s \to \RRR^{m_1 + m_2}$ and $h: \RRR^n \times \RRR^s \to \RRR^{k}$ are mappings, $A \in \RRR^{k \times (m_1+m_2)}$ and $D$ is a square matrix with $(m_1+m_2)$ rows and columns. We will assume $D$ to be positive definite.

\medskip

To our knowledge, stability of two-stage stochastic programs with mixed-integer quadratic recourse has only been studied in \cite{ChenXu2009}. While the authors of the mentioned paper also consider the situation where $D$ is positive semidefinite, they fix a compact set $\Xi \subset \mathbb{R}^{s}$ and confine their stability analysis to the class $\mathcal{P}_{\Xi}$ of Borel probability measures $\mu$ on $\RRR^s$ that satisfy $\mu[\Xi] = 1$. Note that for any gauge function $\psi: \mathbb{R}^{s} \to [0, \infty)$ and any compact set $\Xi \subset \mathbb{R}^{s}$, $\mathcal{P}_{\Xi} \subseteq \mathcal{M}_{1}^{\psi}(\RRR^s)$ is relatively compact for the $\psi-$weak topology (cf. \cite[Corollary A.47]{FoellmerSchied2011}), and hence locally uniformly $\psi-$integrating (see Lemma \ref{relative compactness} below).

\medskip

Furthermore, \cite{ChenXu2009} only examines an objective function that is based on the expectation (although their model may reflect some kind of risk-aversion, see page 465 in \cite{ChenXu2009} for details) and assumes $d$ and $h$ to be of a special form (for any fixed $z$, $d$ is linear in $x$, while $h$ does not depend on $x$). Consequently, the present analysis allows some of the existing results in various directions.

\medskip

\begin{assumption} \label{AMIQP}
\begin{itemize}
\item $A$ and $D$ have rational entries,
\item $D$ is symmetric and positive definite,
\item $\lbrace y \in \RRR^{m_1} \times \mathbb{Z}^{m_2} \mid Ay \leq t \rbrace \neq \emptyset$ for any $t \in \RRR^k$,
\item the mappings $q$ and $h$ are continuous,
\item $\|q(\cdot, \cdot)\|$ and $\|h(\cdot, \cdot)\|$ are limited by exponents $\gamma_q, \gamma_h \in ]0,\infty[$, respectively.
\end{itemize}
\end{assumption}

\medskip

\begin{Prop}
Under Assumption \ref{AMIQP}, the mapping $f$ defined in \eqref{MIQPf} is finite, lower semicontinuous and limited by the exponent $\max \lbrace 2 \gamma_q, 2\gamma_h \rbrace$. In particular, Assumption \ref{Af} is fulfilled. Moreover, $f$ is continuous if
$$
\left\{y_2 \in \mathbb{Z}^{m_2}  \; | \;  \exists y_1 \in \mathbb{R}^{m_1}: \; A\begin{pmatrix} y_1 \\ y_2 \end{pmatrix} \leq h(x,z) \right\}
$$
does not depend on $(x,z)$.
\end{Prop}

\medskip

\begin{proof}
Finiteness and lower semicontinuity of $f$ follow from \cite[Theorem 2.2]{ChenHan2010} and \cite[Lemma 2.7, Remark 2.8]{ChenXu2009}. By the same results, there exists a constants $\alpha, \beta > 0$ such that
\begin{equation*}
|f(x,z) - f(x',z')| \leq \alpha \max\lbrace \|H\|, \|Q\|, \|H'\|, \|Q'\| \rbrace (\|H - H'\| + \|Q - Q'\| + 1) + \beta,
\end{equation*}
where $H= h(x,z)$, $H'= h(x',z')$, $Q= q(x,z)$ and $Q'= q(x',z')$, holds for any $(x,z), (x',z') \in \RRR^n \times \RRR^s$. Consequently, we have
\begin{equation*}
|f(x,z)| \leq \alpha (\|h(x,z)\| + \|q(x,z)\| + 1)^2 + \beta + |f(0,0)|
\end{equation*}
for any $(x,z) \in \RRR^n \times \RRR^s$. Under Assumption \ref{AMIQP}, the latter implies that $f$ is limited by the exponent $\max \lbrace 2 \gamma_q, 2\gamma_h \rbrace$. The second part of the Proposition follows from \cite[Lemma 2.9]{ChenXu2009}.
\end{proof}

\subsection{Mixed-integer convex recourse}

Finally, we consider the fairly general class of mixed-integer problems where the continuous relaxation is convex. Let $f$ be given by
\begin{equation}
\label{Convexf}
f(x,z) = \inf_{y} \lbrace  v(y) \; | \; g(y) \leq h(x,z), y \in \mathbb{R}^{m_1} \times \mathbb{Z}^{m_2} \rbrace,
\end{equation}
where $v: \mathbb{R}^{m_1 + m_2} \to \mathbb{R}$ is convex, the right-hand side of the constraint system is given by the mapping $h: \mathbb{R}^{n} \times \mathbb{R}^{s} \to \mathbb{R}^{k}$ and $g = (g_1, \ldots, g_k)^\top: \mathbb{R}^{m_1 + m_2} \to \mathbb{R}^{k}$ is such that for any $i \in \lbrace 1, \ldots, k \rbrace$, $g_i$ has a closed and convex epigraph.

\medskip

As far as we know, there is no systematic investigation of stability of two-stage stochastic programs where the recourse is given by \eqref{Convexf}. Under assumptions involving a compactness condition, we shall show that $f$ is finite, lower semicontinuous and limited by finite exponent. Hence, Corollary \ref{CorStability} can be applied to derive qualitative stability of the resulting problem \eqref{StochasticProblem} whenever Assumption \ref{Arho} is fulfilled.

\medskip

Let $C(x,z)$ denote the feasible set of the problem in \eqref{Convexf}, i.e. define
$$
C(x,z) := \lbrace  y \in \mathbb{R}^{m_1} \times \mathbb{Z}^{m_2} \; | \;  g(y) \leq h(x,z) \rbrace.
$$
In addition, let us consider the sets 
$$
C_{rel}(t) := \lbrace y \in \RRR^{m_1 + m_2} \mid g(y) \leq t \rbrace, \; \; \; t \in \RRR^k.
$$
Obviously, $C(x,z) = C_{rel}(h(x,z)) \cap (\mathbb{R}^{m_1} \times \mathbb{Z}^{m_2})$ holds for any $(x,z) \in \RRR^n \times \RRR^s$. The following is known about $C_{rel}(\cdot)$.

\medskip

\begin{lemma}[\textbf{\cite[Corollary 5]{AuslenderCrouzeix1988}}] \label{TheoremAuslenderCrouzeix} \ \\
Assume that $C_{rel}(t) \neq \emptyset$ holds for any $t \in \mathbb{R}^{k}$ and that $C_{rel}(0)$ is compact. Let $t_1, \ldots, t_k$ denote the components of $t \in \mathbb{R}^{k}$, then for any $r > 0$,
\begin{equation*}
K(r) := \sup_{t: \|t\| < r, \; y \notin C_{rel}(t)} \frac{\inf \lbrace \| y - y'\| \; | \; y' \in C_{rel}(t) \rbrace}{\max \lbrace g_{j}(y) - t_{j} \; | \; j = 1, \ldots, k \rbrace}
\end{equation*}
is finite and such that 
\begin{equation*}
d_{\infty}(C_{rel}(t), C_{rel}(t')) \leq K(r) \|t - t'\|,
\end{equation*}
holds for any $t,t' \in \RRR^k$ satisfying $ \|t\|, \|t'\| < r$. Here, $d_{\infty}$ denotes the Hausdorff distance.
\end{lemma}

\medskip

We shall work with the following assumptions.

\medskip

\begin{assumption} \label{AConv}
\begin{itemize}
\item The mapping $v$ is convex,
\item for every $i=1, \ldots, k$, the epigraph of $g_i$ is closed and convex,
\item $C_{rel}(t) \cap (\mathbb{R}^{m_1} \times \mathbb{Z}^{m_2}) \neq \emptyset$ for any $t \in \RRR^k$,
\item $C_{rel}(0)$ is compact,
\item the mapping $h$ is continuous,
\item $\|h(\cdot, \cdot)\|$ is limited by an exponent $\gamma_h \in ]0,\infty[$,
\item there are positive constants $\gamma_v, \kappa_v$ such that $|v(y)| \leq \kappa_{v}(\|y\|^{\gamma_{v}} + 1)$ holds for any $y \in \RRR^{m_1+m_2}$,
\item there are positive constants $\gamma_K, \kappa_K$ such that $K(r) \leq \kappa_{K}(r^{\gamma_{K}} + 1)$ holds for any $r > 0$.
\end{itemize}
\end{assumption}

\medskip

The technical conditions in the last part of Assumption \ref{AConv} allow us to prove that $f$ is limited by a finite exponent.

\medskip

\begin{Prop}
Under Assumption \ref{AConv}, the mapping $f$ defined in \eqref{Convexf} is finite, lower semicontinuous and limited by the exponent $\gamma_{h}(\gamma_{K} + 1)(\gamma_{v} + 1)$. In particular, Assumption \ref{Af} is fulfilled. Moreover, $f$ is continuous if $C(x,z) = C_{rel}(h(x,z))$ holds for any $(x,z) \in \RRR^n \times \RRR^s$.
\end{Prop}

\medskip

\begin{proof} Fix any $(x,z) \in \mathbb{R}^{n} \times \mathbb{R}^{s}$. Since $g$ is convex and hence continuous, $C_{rel}(h(x,z))$ is closed and the boundedness of $C_{rel}(0)$ and Theorem \ref{TheoremAuslenderCrouzeix} yield that $C_{rel}(h(x,z))$ is bounded. Consequently, $C(x,z)$ is the intersection of a compact and closed set and thus compact. Furthermore, $v$ is convex and hence continuous, which implies that
\begin{equation*}
\inf_{y} \lbrace v(y) \; | \; y \in C(x,z) \rbrace
\end{equation*}
is finite and that the infimum is attained. Hence, $f$ is real-valued and admits the representation $f(x,z) = \min_{y} \lbrace v(y) \; | \; y \in C(x,z) \rbrace$.

\medskip

Next, we shall prove that the set-valued mapping $C: \RRR^n \times \RRR^s \to 2^{\RRR^{m_1} \times \mathbb{Z}^{m_2}}$ is upper semicontinuous: Otherwise, there would exist a point $(x_{0},z_{0}) \in \mathbb{R}^{n} \times \mathbb{R}^{s}$, an open set $\mathcal{O} \subseteq \mathbb{R}^{m_1+m_2}$ and sequences $\lbrace (x_{l}, z_{l}) \rbrace_{l \in \mathbb{N}} \subseteq \mathbb{R}^{n} \times \mathbb{R}^{s}$ and $\lbrace (y_{l}) \rbrace_{l \in \mathbb{N}} \subseteq \mathbb{R}^{m_1+m_2}$ such that $C(x_{0},z_{0}) \subset \mathcal{O}$, $y_{l} \in C(x_{l}, z_{l})$, $y_{l} \notin \mathcal{O}$,
$$
\|(x_l,z_l)-(x_0, z_0)\| \leq \frac{1}{l} \; \text{and} \; \|h(x_l,z_l)-h(x_0, z_0)\| \leq \frac{1}{l} 
$$
hold for any $l \in \mathbb{N}$. By Theorem \ref{TheoremAuslenderCrouzeix} we have
\begin{equation*}
\sup_{l \in \mathbb{N}} d_{\infty}(C_{rel}(h(x_{0}, z_{0})), C_{rel}(h(x_{l}, z_{l}))) \leq K(\|h(x_0,z_0)\| + 1) =: K_0,
\end{equation*}
which implies that
\begin{equation*}
\bigcup_{l = 1}^{\infty} \lbrace y_{l} \rbrace \; \subseteq \; \bigcup_{l = 1}^{\infty} C(x_{l},z_{l}) \; \subseteq \; \bigcup_{l = 1}^{\infty} C_{rel}(h(x_{l},z_{l})) \; \subseteq \; \lbrace C_{rel}(h(x_{0},z_{0})) \rbrace + B_{K_0}(0)
\end{equation*}
is bounded. Here, $B_{K_0}(0) \subset \RRR^{m_1 + m_2}$ denotes the open $\|\cdot\|-$ball of radius $K_0$ centered at $0$. Consequently, we can assume $y_{l} \to \bar{y}$ for some $\bar{y} \in \mathbb{R}^{m_1 + m_2}$ without loss of generality. $y_{l} \in C(x_{l},z_{l}) \subseteq \mathbb{R}^{m_{1}} \times \mathbb{Z}^{m_{2}}$ holds for any $l \in \mathbb{N}$ and implies $\bar{y} \in \mathbb{R}^{m_{1}} \times \mathbb{Z}^{m_{2}}.$ Furthermore, by $h(x_{l},z_{l}) \to h(x_{0},z_{0})$, Theorem \ref{TheoremAuslenderCrouzeix} yields that $\bar{y} \in C_{rel}(h(x_{0},z_{0}))$. Thus,
\begin{equation*}
\bar{y} \in C_{rel}(h(x_{0},z_{0})) \cap (\mathbb{R}^{m_{1}} \times \mathbb{Z}^{m_{2}}) = C(x_{0},z_{0}) \subset \mathcal{O}.
\end{equation*}
On the other hand, $y_{l} \notin \mathcal{O}$ for any $l \in \mathbb{N}$ and $\mathcal{O}$ is open, which yields the contradiction $\bar{y} \notin \mathcal{O}$. Hence, $C$ is upper semicontinuous. Since $h$ and $v$ are continuous, the upper semicontinuity of $C$ allows us to apply \cite[Theorem 2]{Berge1963} to conclude that $f$ is lower semicontinuous and hence Borel measurable.

\medskip

Next, we shall prove that $f$ is limited by the exponent $\gamma_{h}(\gamma_{K} + 1)(\gamma_{v} + 1)$: Let $(x,z) \in \mathbb{R}^{n} \times \mathbb{R}^{s}$ be fixed. Based on the considerations above, there exist $y^{*} \in C(x,z)$ and $y_{0}^{*} \in C(0,0)$ satisfying $f(x,z) = v(y^{*})$ and $f(0,0) = v(y_{0}^{*})$. Since $C_{rel}(h(0,0))$ is compact,
\begin{equation*}
d_{0} := \max \lbrace \|y - y'\| \; | \; y, y' \in C_{rel}(h(0,0)) \rbrace
\end{equation*}
is finite. Theorem \ref{TheoremAuslenderCrouzeix} and Assumption \ref{AConv} yield
\begin{align*}
\|y^{*} - y^{*}_{0}\| \; \leq \; &d_{\infty}(C_{rel}(h(x,z)),C_{rel}(h(0,0))) + d_{0} \\
\leq \; &K(\|h(x,z)\| + \|h(0,0)\| + 1)(\|h(x,z)\| + \|h(0,0)\|) + d_{0} \\
\leq \; &2 \kappa_{K}(\|h(x,z)\| + \|h(0,0)\| + 1)^{\gamma_{K} + 1} + d_{0} \\
\leq \; &\kappa^{*}(\|h(x,z)\|^{\gamma_{K} + 1} + 1),
\end{align*}
where $\kappa^{*} := 2^{\gamma_{K} + 2} \kappa_{K} (\|h(0,0)\| +1)^{\gamma_{K} + 1} + d_{0}$. For any $r>0$, the convex function $v$ is Lipschitz continuous on the open $\|\cdot\|-$ball of radius $r$ centered at $0$ by \cite[Theorem A, Lemma A]{RobertsVarberg1974} and the Lipschitz constant is given by
\begin{equation*}
L_{v}(r) := \frac{2}{r} \left( \max_{y \in \lbrace 2r, -2r \rbrace^{m_1+m_2}} |v(y)| + 2|v(0)| \right).
\end{equation*}
Using Assumption \ref{AConv}, we obtain
\begin{align*}
L_{v}(\|y^{*} - y^{*}_{0}\| + \|y^{*}_{0}\| + 1) \; \leq \; &\max_{y \in \lbrace \pm 2(\|y^{*} - y^{*}_{0}\| + \|y^{*}_{0}\| + 1) \rbrace^{m_1+m_2}} 2|v(y)| + 4|v(0)| \\
\leq \; &4\kappa_{v}(2\sqrt{m_1+m_2})^{\gamma_{v}}(\|y^{*} - y^{*}_{0}\| + \|y^{*}_{0}\| + 1)^{\gamma_{v}} + 4|v(0)| \\
\leq \; &\kappa_{L}(\|y^{*} - y^{*}_{0}\|^{\gamma_{v}} + 1),
\end{align*}
where $\kappa_{L} := 4\kappa_{v}(4\sqrt{m_1+m_2})^{\gamma_{v}}(\|y^{*}_{0}\| + 1)^{\gamma_{v}} + 4|v(0)|$. Since 
\begin{equation*}
\|y^{*}\| < \|y^{*} - y^{*}_{0}\| + \|y^{*}_{0}\| + 1,
\end{equation*}
the above considerations yield
\begin{align*}
|f(x,z)| \; \leq \; &|v(y^{*}) - v(y^{*}_{0})| + |v(y^{*}_{0})| \\
\leq \; &L_{v}(\|y^{*} - y^{*}_{0}\| + \|y^{*}_{0}\| + 1)\|y^{*} - y^{*}_{0}\| + |v(y^{*}_{0})| \\
\leq \; &2 \kappa_{L} \|y^{*} - y^{*}_{0}\|^{\gamma_{v} + 1} + 2 \kappa_{L} + |v(y^{*}_{0})| \\
\leq \; &2 \kappa_{L} (\kappa^{*})^{\gamma_{v} + 1} (\|h(x,z)\|^{\gamma_{K} + 1} + 1)^{\gamma_{v} + 1} + 2 \kappa_{L} + |v(y^{*}_{0})| \\
\leq \; &\bar{\kappa}(\|h(x,z)\|^{(\gamma_{K} + 1)(\gamma_{v} + 1)} + 1),
\end{align*}
where $\bar{\kappa}:= 2^{\gamma_{v} + 2} \kappa_{L} (\kappa^{*})^{\gamma_{v} + 1} + 2 \kappa_{L} + |v(y^{*}_{0})|$. Invoking that $\|h(\cdot,\cdot)\|$ is limited by $\gamma_h$, we finally may conclude that $f$ is limited by the exponent $\gamma_{h}(\gamma_{K} + 1)(\gamma_{v} + 1)$.

\medskip

The continuity of $f$ in the second part of the proposition follows from \cite[Theorem 2]{Berge1963} and the fact that $C = C_{rel} \circ h$ is both upper and lower semicontinuous by Theorem \ref{TheoremAuslenderCrouzeix}.
\end{proof}

\medskip

\section{$\psi-$weak topology and the topology of weak convergence}
\label{psi-weak topology}

Throughout this section we want to gather some useful results on the $\psi$-weak topologies we have already introduced in Section \ref{mainresult}. So let us fix a gauge function $\psi: \RRR^{d}\rightarrow [0,\infty[$, and recall that $\cM_{1}^{\psi}(\RRR^{d})$ denotes the set of all Borel probability measures $\mu$ on $\RRR^{d}$ satisfying $\int\psi\,d\mu<\infty$. The $\psi-$weak topology on $\cM_{1}^{\psi}(\RRR^{d})$ may be characterized in the following ways.
\medskip

\begin{lemma}
\label{psischwacheTopologie}
The $\psi$-weak topology is metrizable, and for every sequence $(\mu_{n})_{n\in\NNN_{0}}$ in $\cM_{1}^{\psi}(\RRR^{d})$ the following statements are equivalent.
\begin{enumerate}
\item [(1)] $\mu_{n}\to\mu_{0}$ w.r.t. the $\psi-$weak topology.
\item [(2)] $\mu_n\to\mu_{0}$ w.r.t. the topology of weak convergence, and $\int_{\RRR^{d}}\psi~d\mu_{n}\to \int_{\RRR^{d}}\psi~d\mu_{0}.$
\item [(3)] $\mu_n\to\mu_{0}$ w.r.t. the topology of weak convergence, and 
$$
\lim_{a\to\infty}\sup_{n\in\NNN}\int_{\RRR^{d}}\psi\cdot\eins_{]a,\infty[}(\psi)~d\mu_{n}= 0.
$$
\end{enumerate}
In particular, the $\psi$-weak topology is metrizable by 
$$
d_{\psi}: \cM_{1}^{\psi}(\RRR^{d})\times\cM_{1}^{\psi}(\RRR^{d})\rightarrow\RRR,~(\mu,\nu)\mapsto d(\mu,\nu) + \left|\int_{\RRR^d} \psi~d\mu - \int_{\RRR^d} \psi~d\nu\right|,
$$
where $d$ denotes any metric which generates the topology of weak convergence, e.g. the Prokhorov metric.
\end{lemma}
\medskip

\begin{proof}
From Theorem A.38 in \cite{FoellmerSchied2011} it is known that the the $\psi-$weak topology is metrizable. The equivalence of (1) and (2) has been shown in
\cite[Lemma 3.4]{KraetschmerEtAll2012}. Finally, the equivalence of (2) and (3) follows immediately from the convergence of moments theorem (cf. \cite[Theorem 2.20]{vanderVaart1998}). Finally, $d_{\psi}$ obviously defines a metric, which by (2) generates the $\psi$-weak topology.
\end{proof}

\medskip

We are interested in gauge functions of the form $\psi := \|\cdot\|^{q}$ for any $q > 0,$ where $\|\cdot\|$ denotes the Euclidean norm on $\RRR^{d}$. For $q\geq 1,$ the $\|\cdot\|^{q}-$weak topology is generated by so called \textit{Wasserstein metric $d_{W,d,q}$ of order $q$} w.r.t. $\|\cdot\|^{q}$ defined by
\begin{eqnarray*}
d_{W,d,q}(\mu,\nu) 
&:=& 
\inf\left\{\left(\int_{\RRR^{d}\times\RRR^{d}}\|x - y\|^{q}~\pi(dx,dy)\right)^{1/q}\mid \pi\in \cM_{1}(\mu,\nu)\right\}\\[0.2cm]
&=&
\left(\inf\left\{\int_{\RRR^{d}\times\RRR^{d}}\|x - y\|^{q}~\pi(dx,dy)\mid \pi\in \cM_{1}(\mu,\nu)\right\}\right)^{1/q},
\end{eqnarray*}
where $\cM_{1}(\mu,\nu)$ denotes the set of all Borel probability measures on $\RRR^{d}\times\RRR^{d}$ with $\mu$ as the first $d-$dimensional marginal and $\nu$ as the second one (cf. \cite[Theorem 6.3.1]{Rachev1991}). 
The $\|\cdot\|^{q}-$weak topology may be also generated by the {\textit{Fortet-Mourier metric $d_{FM, d, q}$ of order q}} w.r.t. $\|\cdot\|^{q}$ defined by
\begin{eqnarray*}
&&
d_{FM,d,q}(\mu,\nu)\\ 
&:=&
\inf\left\{\int_{\RRR^{d}\times\RRR^{d}}\|x - y\|~\max\{1, \|x\|^{q-1},\|y\|^{q-1}\}~\pi(dx,dy)\mid \pi\in \cL_{1}(\mu,\nu)\right\},
\end{eqnarray*}
where $\cL_{1}(\mu,\nu)$ denotes the set of all finite Borel measures on $\RRR^{d}\times\RRR^{d}$ satisfying 
$$\pi(A\times\RRR^{d}) - \pi(\RRR^{d}\times A) = \mu(A) - \nu(A)\quad\mbox{for every Borel measurable set}~ A\subseteq\RRR^{d}
$$ 
(see \cite[Theorem 6.3.1]{Rachev1991}).
\medskip

Next we want to identify those subsets $\cM\subseteq\cM_{1}^{\|\cdot\|^{q}}(\RRR^{d})$ on which the $\|\cdot\|^{q}-$ weak topology and the topology of weak convergence coincide. Analogously to Section \ref{mainresult} we shall call a subset $\cM\subseteq\cM_{1}^{\|\cdot\|^{q}}(\RRR^{d})$ {\textit{locally uniformly $\|\cdot\|^{q}-$ integrating}} if for any $\nu\in\cM$ and every $\varepsilon > 0$, there exists some open neighbourhood $\cL$ of $\nu$ w.r.t. the topology of weak convergence such that
$$
\lim_{a\to\infty}\sup_{\mu\in\cL\cap\cM}\int_{\RRR^{d}}~\|z\|^{q}\cdot\eins_{]a,\infty[}(\|z\|^{q})~\mu(dz)\leq\varepsilon.
$$ 
The following result may be found in \cite{Zaehle2016} (Lemma 3.4 there).
\begin{proposition}
\label{coincidence}
For $q > 0$ and $\cM\subseteq\cM_{1}^{\|\cdot\|^{q}}(\RRR^{d})$ the following statements are equivalent.
\begin{enumerate}
\item [(1)] The $\|\cdot\|^{q}-$weak topology and the topology of weak convergence coincide on $\cM$. 
\item [(2)] $\cM$ is locally uniformly $\|\cdot\|^{q}-$integrating.
\end{enumerate}
\end{proposition}
\medskip

Next, we shall provide a useful criterion to verify a set of Borel probability measures on $\RRR^{d}$ as locally uniformly $\|\cdot\|^{q}-$ integrating.
\begin{lemma}
\label{relative compactness}
A subset $\cM\subseteq\cM_{1}^{\|\cdot\|^{q}}(\RRR^{d})$ that is relatively compact for the $\|\cdot\|^{q}$-weak topology is locally uniformly $\|\cdot\|^{q}-$integrating.
\end{lemma}
\medskip

\begin{proof}
The statement of Lemma \ref{relative compactness} is a special case of Lemma 3.1 from 
\cite{KraetschmerEtAll2015b}.
\end{proof}

\section{Proof of Theorem \ref{QRestriktionStetig}}

Let $\rho$ be as in Assumption \ref{Arho} and consider the mapping 
$$
\Theta_{f}: \RRR^{n}\times \mathcal{M}^{\| \cdot \|^{\gamma p}}_{1}(\mathbb{R}^{s}) \rightarrow \mathcal{M}^{| \cdot |^{p}}_{1}(\mathbb{R}),~(x,\nu)\mapsto  (\delta_{x} \otimes \nu) \circ f^{-1}
$$

which is well-defined under Assumption \ref{Af} due to Lemma \ref{LemmaEndlichkeit}. Let us equip $\cM^{\|\cdot\|^{\gamma p}}(\RRR^{s})$ and $\cM^{|\cdot|^{p}}(\RRR)$ respectively with the $\|\cdot\|^{\gamma p}-$weak topology $\tau_{s,\gamma p}$ and the $|\cdot|^{p}-$weak topology $\tau_{1,p}$ as defined in Section \ref{psi-weak topology}. Furthermore, let $\tau_{\RRR^{n}}\otimes\tau_{s,\gamma p}$ denote the product topology of the standard topology on $\RRR^{n}$ and $\tau_{s,\gamma p}.$ An important step in our argumentation is to investigate continuity of the mapping $\Theta_{f}$ w.r.t. $\tau_{\RRR^{n}}\otimes\tau_{s,\gamma p}$ and $\tau_{1,p}$.

\medskip

\begin{proposition} 
\label{PhiFStetig}
Let Assumption \ref{Af} be fulfilled and let $x\in\RRR^{n}$ and $\nu\in\mathcal{M}^{\| \cdot \|^{\gamma p}}_{1}(\mathbb{R}^{s})$ be such that $\delta_{x} \otimes \nu(D_{f}) = 0$, where $D_{f}$ is the set of discontinuity points of $f$. Then $\Theta_{f}$ is continuous at $(x,\nu)$ w.r.t. $\tau_{\RRR^{n}}\otimes\tau_{s,\gamma p}$ and $\tau_{1,p}$.
\end{proposition}

\medskip

\begin{proof}
Since $\tau_{\RRR^{n}}\otimes\tau_{s,\gamma p}$ and $\tau_{1,p}$ are metrizable, it suffices to show that $\Phi_{f}$ is sequentially continuous at $(x,\nu)$. So let $(x_l)_{l\in\NNN}$ be any sequence in $\RRR^{n}$ converging to $x$, and let $(\nu_{l})_{l \in \mathbb{N}}$ be a sequence in $\mathcal{M}^{\| \cdot \|^{\gamma p}}_{1}(\mathbb{R}^{s})$ such that $\nu_{l}\to\nu$ $\| \cdot \|^{\gamma p}-$weakly.

\medskip

Then obviously $\delta_{x_l} \otimes \nu_l \to \delta_x \otimes \nu$ with respect to the topology of weak convergence on $\RRR^{n}\times\RRR^{s}.$ By assumption, $\delta_{x} \otimes \nu(D_{f}) = 0$ holds so that we may conclude from the continuous mapping theorem (see e.g. \cite[Theorem 5.1]{Billingsley1968}) that $(\delta_{x_l} \otimes \nu_l) \circ f^{-1} \to (\delta_x \otimes \nu) \circ f^{-1}$ w.r.t. the topology of weak convergence on $\RRR.$ 
\medskip

Furthermore, as in the proof of Lemma \ref{LemmaEndlichkeit} we obtain
$$
|f(x_{l},z)|^{p}\leq \eta(x_{l})^{p} 2^{p}(\|z\|^{p\gamma} + 1)\quad\mbox{for}~l\in\NNN,~z\in\RRR^{s}
$$
and some locally bounded mapping $\eta$. Without loss of generality we may assume that $C:= \sup_{l\in\NNN}\eta(x_{l}) < \infty.$ Thus by Fubini-Tonelli theorem
\begin{eqnarray*}
&&
\sup_{l\in\NNN}\int_{\mathbb{R}^{n} \otimes \mathbb{R}^{s}} 
|f|^{p}\eins_{]a,\infty[}(|f|^{p})~d(\delta_{x_l} \otimes \nu_l)\\ 
&\leq&
 2^{p}C^{p}\sup_{l\in\NNN}\int_{\mathbb{R}^{n} \times \mathbb{R}^{s}} 
(\|z\|^{\gamma p} + 1)\eins_{]a/(2^{p} C^{p}),\infty[}(\|z\|^{\gamma p} + 1)~(\delta_{x_l} \otimes \nu_l)(dx,dz)\\
&=&
2^{p}C^{p} \sup_{l\in\NNN}\int_{\mathbb{R}^{s}} 
(\|z\|^{\gamma p} + 1)\eins_{]a/(2^{p} C^{p}),\infty[}(\|z\|^{\gamma p} + 1)~\nu_l(dz)
\end{eqnarray*}
Since $\nu_{l}\to\nu$ $\|\cdot\|^{\gamma p}-$weakly, we may conclude from Lemma \ref{psischwacheTopologie}
$$
\lim_{a\to\infty}\sup_{l\in\NNN}\int_{\mathbb{R}^{s}} 
(\|z\|^{\gamma p} + 1)\eins_{]a/(2^{p} C^{p}),\infty[}(\|z\|^{\gamma p} + 1)~ \nu_l(dz) = 0.
$$
This implies
$$
\lim_{a\to\infty}\sup_{l\in\NNN}\int_{\mathbb{R}} |\cdot|^{p}~d(\delta_{x_l} \otimes \nu_l) \circ f^{-1}
=
\lim_{a\to\infty}\sup_{l\in\NNN}\int_{\mathbb{R}^{n} \otimes \mathbb{R}^{s}} 
|f|^{p}\eins_{]a,\infty[}(|f|^{p})~d(\delta_{x_l} \otimes \nu_l) = 0,
$$
so that by Lemma \ref{psischwacheTopologie} again
\begin{eqnarray*}
\lim_{l\to\infty}\int_{\mathbb{R}} |\cdot|^{p}~d(\delta_{x_l} \otimes \nu_l) \circ f^{-1} 
&=& 
\lim_{l\to\infty}\int_{\mathbb{R}^{n} \times \mathbb{R}^{s}} |f|^{p}~d(\delta_{x_l} \otimes \nu_l)\\ 
&=& 
\int_{\mathbb{R}^{n} \times \mathbb{R}^{s}} |f|^{p}~d(\delta_{x} \otimes \nu)
=
\int_{\mathbb{R}} |\cdot|^{p}~d(\delta_{x} \otimes \nu) \circ f^{-1}.
\end{eqnarray*}

Therefore $(\delta_{x_l} \otimes \nu_l) \circ f^{-1} \to (\delta_x \otimes \nu) \circ f^{-1}$ w.r.t. the $|\cdot|^{p}-$weak topology.
\end{proof}
\medskip

The line of reasoning to prove Theorem \ref{QRestriktionStetig} may be described roughly as follows. We want to verify continuity of 
$$
Q:\RRR^{n}\times \mathcal{M}^{\| \cdot \|^{\gamma p}}_{1}(\mathbb{R}^{s}) \to \mathbb{R},~(x,\nu)\mapsto {\cal R}_{\rho}\left((\delta_{x}\otimes\nu)\circ f^{-1}\right) = ({\cal R}_{\rho} \circ \Theta)(x,\nu)
$$
w.r.t. the production topology of the standard topology on $\RRR^{n}$ and the $\|\cdot\|^{\gamma p}-$weak topology on $\cM_{1}^{\|\cdot\|^{\gamma p}}(\RRR^{s})$, invoking then Proposition \ref{coincidence}. In view of Proposition \ref{PhiFStetig} the first part may be concluded if the continuity of $\cR_{\rho}$ w.r.t. the $|\cdot|^{p}$-weak topology can be shown. This will be done next.
\medskip

\begin{theorem} \label{RRhoStetig}
Under Assumption \ref{Arho}, $\cR_{\rho}$ is continuous with respect to the $|\cdot|^{p}$-weak topology, and it is nondecreasing w.r.t. the increasing convex order if in addition $\rho$ is translation-equivariant.
\end{theorem}
\medskip

\begin{proof}
Let us first recall that $L^{p}(\Omega, \mathcal{F}, \mathbb{P})$, equipped with ordinary $L^{p}-$norm $\|\cdot\|_{p}$ and the $\mathbb{P}$-almost sure partial order, is a Banach lattice, i.e. $\|\cdot\|_{p}$ is complete satisfying $\|X\|_{p}\leq \|Y\|_{p}$ whenever 
$|X|\leq |Y|$ $\mathbb{P}-$a.s. (see Theorem 13.5 in \cite{AliprantisBorder2006}). Furthermore, $\rho$ is assumed to be a convex mapping on the Banach lattice $L^{p}\OFP$ which is nondecreasing w.r.t. the $\mathbb{P}$-almost sure partial order. Hence it is continuous w.r.t. $\|\cdot\|_{p}$ (cf. \cite[Corollary 2]{BiaginiFritelli2009}).

\medskip

Since the $|\cdot|^{p}-$weak topology is metrizable it suffices to show sequential continuity for ${\cal R}_{\rho}.$ So consider any sequence $(\mu_l)_{l \in \mathbb{N}}$ in $\mathcal{M}^{| \cdot |^{p}}_{1}(\mathbb{R})$ which converges to some $\mu \in \mathcal{M}^{| \cdot |^{p}}_{1}(\mathbb{R})$ w.r.t. the $|\cdot|^{p}$-weak topology. Then by Theorem 3.5 in \cite{KraetschmerEtAll2014} we may  find a sequence $(X_l)_{l \in \mathbb{N}}$ in $L^{p}(\Omega, \mathcal{F}, \mathbb{P})$ and some $X \in L^{p}(\Omega, \mathcal{F}, \mathbb{P})$ such that $\mu_l$ is the distribution of $X_l$ for $l\in\NNN,$  $X$ has $\mu$ as its distribution, and $\int_{\Omega} |X_l - X|^{p}~d\mathbb{P} \to 0$. 
So by law-invariance and $L^{p}$-norm-continuity of $\rho$ we have
\begin{equation*}
\cR_{\rho}(\mu_l) = \rho(X_l) \to \rho(X) = \cR_{\rho}(\mu). 
\end{equation*}
This shows continuity of $\cR_{\rho}$ w.r.t. the $|\cdot|^{p}$-topology. For the second part of the statement let us assume that $\rho$ is translation-equivariant. Then due to $L^{p}$-norm continuity of $\rho$ the application of \cite[Corollary 4.65]{FoellmerSchied2011} yields that the restriction of $\cR_{\rho}$ to the set of probability measures on $\RRR$ with bounded support is nondecreasing w.r.t. the increasing convex order. Then invoking $L^{p}$-norm continuity of $\rho$ again, we obtain that $\cR_{\rho}$ is nondecreasing w.r.t. the increasing convex order on its entire domain. This completes the proof.
\end{proof}

\medskip

Combining Proposition \ref{PhiFStetig} with Theorem \ref{RRhoStetig}, we obtain immediately the following criterion to guarantee continuity of the function $Q$ w.r.t. the production topology of the standard topology on $\RRR^{n}$ and the $\|\cdot\|^{\gamma p}-$weak topology on $\cM_{1}^{\|\cdot\|^{\gamma p}}(\RRR^{s}).$

\medskip

\begin{theorem} \label{QStetig}
Let Assumptions \ref{Arho} and \ref{Af} be fulfilled. If $x \in \mathbb{R}^{n}$ and $\nu \in \cM_{1}^{\|\cdot\|^{\gamma p}}(\RRR^{s})$ satisfy $\delta_{x} \otimes \nu(D_{f}) = 0$, the mapping $Q$ is continuous at $(x, \nu)$ with respect to the product topology of the standard topology on $\mathbb{R}^{n}$ and the $\|\cdot\|^{\gamma p}-$weak topology on $\cM_{1}^{\|\cdot\|^{\gamma p}}(\RRR^{s}).$
\end{theorem}

\medskip

Now we are ready to prove the main result.

\medskip

{\textbf{ Proof of Theorem \ref{QRestriktionStetig}:}}\\
Since $\cM$ is assumed to be locally uniformly $\|\cdot\|^{\gamma p}-$integrating, the topology of weak convergence and the $\|\cdot\|^{\gamma p}$-weak topology coincide on $\cM$ due to Proposition \ref{coincidence}. Thus Theorem \ref{QRestriktionStetig} follows immediately from Theorem \ref{QStetig}.
\vbox{\hrule height0.6pt\hbox{%
   \vrule height1.3ex width0.6pt\hskip0.8ex
   \vrule width0.6pt}\hrule height0.6pt
  }

\section{Concluding Remarks}
We offer a general framework to derive stability results for mean-risk models with respect to pertubations of the underlying probability distribution in terms of the topology of weak convergence. Besides unifying already existing stability results, our framework also allows to identify some strategic points to tackle this issue for risk-averse stochastic programming in the following general form: 
$$
\min_{x}~\left\{\cR((\delta_{x}\otimes \nu_{Z})\circ f^{-1})\mid x\in X\right\}.
$$
where 
\medskip

\begin{itemize}
\item $X\subseteq\RRR^{n}$ is nonvoid and compact,\\[-0.2cm]
\item $f:\RRR^{n}\times\RRR^{s}\rightarrow\RRR$ is Borel-measurable and satisfies some polynomial growth condition in its second variable, \\[-0.2cm]
\item $\delta_{x}\otimes\nu_{Z}$ stands for the product measure of the Dirac measure at $x$ and the distribution of a known $s$-dimensional random vector $Z$,\\[-0.2cm]
\item $\cR$ denotes a functional on the set of probability distributions on $\RRR$ with absolute moments of order $p\in [1,\infty[$ representing a law-invariant convex function on an $L^{p}$-space which is monotone nondecreasing w.r.t. the a.s. partial order.
\end{itemize}
\medskip

This formulation enables us to boil down stability to continuity of the mapping
$$
Q:\RRR^{n}\times \mathcal{M} \to \mathbb{R},~(x,\nu)\mapsto  {\cal R}\left((\delta_{x}\otimes\nu)\circ f^{-1}\right)
$$
w.r.t. the product topology of the standard topology on $\RRR^{n}$ and the topology of weak convergence on a set $\mathcal{M}$ containing $\nu_{Z}$, and consisting of probability distributions on $\RRR^{s}$ which fullfill some moment condition corresponding to the growth condition of $f$. The key is that we always have continuity w.r.t. the product topology of the standard topology on $\RRR^{n}$ and a specific topology on $\mathcal{M}$ finer than the topology of weak convergence. As this topology on $\mathcal{M}$ belongs to the general class of the so called $\psi$-weak topologies, we may use a recently obtained concept which describes the sets where $\psi$-weak topologies and the topology of weak convergence coincide. This provides a fairly general condition on ${\cal M}$ which guarantees the desired continuity property, encompassing usually imposed conditions like restrictions on higher moments or fixing a common compact supports for the probability measures from $\mathcal{M}$.
\medskip

The argumentation relies on a specific continuity property of the functional $\mathcal{R}$ which is known to hold also for functionals $\mathcal{R}$ representing law-invariant, almost surely pointwise nondecreasing convex functions on Banach spaces more general than the considered $L^{p}$-spaces (cf. \cite{KraetschmerEtAll2014}, \cite{KraetschmerEtAll2015b}). This offers the opportunity to extend the framework by relaxing the polynomial growth condition on the function $f$ e.g. by means of a general Young function $\Psi$. Then the set $\mathcal{M}$ has to be endowed with a possibly different $\psi$-weak topology. Again we might impose the concept from the theory of $\psi$-weak topologies to identify such sets $\mathcal{M}$ where the chosen $\psi$-weak topology coincides with the topology of weak convergence. However, it is beyond the scope of this paper to work out this sketch of generalization.

\bigskip

{\textbf{Acknowledgements:}} The authors would like to thank two anonymous referees for careful reading, and suggestions which have helped to improve an earlier draft of the paper. The second author is also grateful to Henryk Z\"ahle for fruitful discussions. The first and third authors acknowledge support of the Deutsche Forschungsgemeinschaft within the Research Training Group GRK 1855.
%\textit{Discrete Optimization of Technical Systems under Uncertainty}.
Furthermore, the Deutsche Forschungsgemeinschaft has supported the third author via the Collaborative Research Center TRR 154.
%\textit{Mathematical Modelling, Simulation and Optimization Using the Example of Gas Networks}.

\end{document}